\documentclass[centertags]{article}

\usepackage{amsmath}
\usepackage{amsthm}
\usepackage{amscd}
\usepackage{amssymb}
\usepackage{verbatim}

\title{Group actions on one-manifolds, II: Extensions
of H\"{o}lder's Theorem}

\author{Benson Farb\thanks{Supported in 
part by the NSF and by the Sloan Foundation.}\ \  and 
John Franks\thanks{Supported in part by NSF grant DMS9803346.}}
\date{April 1, 2001}

\newtheorem{theorem}{Theorem}[section]
\newtheorem{proposition}[theorem]{Proposition}
\newtheorem{lemma}[theorem]{Lemma}
\newtheorem{corollary}[theorem]{Corollary}

\newtheorem{defn}[theorem]{Definition}

\newtheorem{exmple}[theorem]{Example}
\newtheorem{thm}[theorem]{Theorem}

\newtheorem{lem}[theorem]{Lemma} 
\newtheorem{prop}[theorem]{Proposition}

\newcommand{\R}{\bf R}
\newcommand{\Z}{{\bf Z}}
\newcommand{\T}{S}

\newcommand{\I}{\mathcal I}

\def\endproof{$\diamond$ \bigskip}
\DeclareMathOperator{\Fix}{Fix}

\DeclareMathOperator{\Homeo}{Homeo}
\DeclareMathOperator{\Aff}{Aff}
\DeclareMathOperator{\Diff}{Diff}
\DeclareMathOperator{\Dist}{Dist}

\DeclareMathOperator{\PSL}{PSL}

\begin{document}
\maketitle
\begin{abstract}
This self-contained paper is part of a series \cite{FF1} seeking
to understand groups of homeomorphisms of manifolds in analogy with
the theory of Lie groups and their discrete subgroups.  In this paper
we consider groups which act on $\R$ with restrictions on the fixed
point set of each element.  One result is a topological
characterization of affine groups in $\Diff^2(\R)$ as those groups whose 
elements have at most one fixed point.
\end{abstract}

\section{Introduction}
In studying finite group actions on manifolds, one typically
relates the algebraic structure of the group, the topology of the
manifold, and the topology of the fixed-set of the action.  An
example of this is Smith theory.  In this paper we try to
understand, in a very simple case, such relationships for actions
of infinite groups.  In particular we extend one of the first
theorems in this direction: H\"{o}lder's theorem.

Let $G$ be a group endowed with an order relation $<$ which is both left
and right invariant.  It is an old theorem of H\"{o}lder \cite{Ho} that
if $(G,<)$ is {\em Archimedean}, that is for every $g,h\in G$ there
exists $n\in \Z$ with $g^n>h$, then $G$ must in fact be abelian.  It was
observed early on that H\"{o}lder's Theorem implies that any group of
homeomorphisms of $\R$ which acts freely must be abelian (see, e.g. 
\cite{Gh} or \cite{FS} for a proof).  This theorem may be extended to
the following.  

We denote by $\Homeo(\R)$ the group of homeomorphisms of the real line
$\R$, and by $\Homeo^+(\R)$ those preserving orientation.  

\begin{theorem}
\label{theorem:compact}
Let $G<\Homeo(\R)$.  Suppose that there is 
a bounded interval $[c,d]$ so that $\Fix(g)\subseteq [c,d]$ for every 
nontrivial $g\in G$.  Then $G$ is abelian.
\end{theorem}

It is not possible to weaken the hypothesis in Theorem
\ref{theorem:compact} by allowing the interval $[c,d]$ to depend on 
the group element $g$: the affine group $\Gamma$ generated by $a(x)=x+1$ 
and $b(x)=2x$ has the property that each element has at most
one fixed point, yet $\Gamma$ is not abelian.

Another variation is to restrict the number but not the location of 
fixed points.

\begin{theorem}
\label{theorem:onefixed}
Let $G<\Homeo(\R)$.  If every nontrivial
element of $G$ has precisely one fixed point, then $G$ is abelian.
\end{theorem}

In fact in this case $G$ is conjugate in $\Homeo(\R)$ to a group of
homeomorphisms fixing the origin and acting on each side of the origin
by multiplication.

\bigskip
\noindent
{\bf Characterizing affine actions. }If one allows every nontrivial 
element of $G$
to have {\em at most} one fixed point, then a much richer collection of
examples appears.  For example the affine group
$\Gamma=<a,b:bab^{-1}=a^2>$ mentioned above acts on $\R$ with every
element having at most one fixed point, yet $\Gamma$ is not abelian.  On
the other hand this group is solvable of derived length two, i.e.\ its
commutator subgroup is abelian.

\begin{theorem}[T. Barbot, N. Kovacevic]
\label{thm:affine}
Let $G<\Homeo^+(\R)$. Suppose that every nontrivial element of $G$ has
at most one fixed point.  Then $G$ is metabelian and, in fact,
isomorphic to a subgroup of the affine group $\Aff(\R)$.
\end{theorem}

Theorem \ref{thm:affine} provides an abstract isomorphism between $G$
and an affine group but does not relate the action to the standard
action of $\Aff(\R)$ on $\R.$ However, using a quasi-invariant measure
technique of Plante \cite{P} it is not difficult to prove the following.

\begin{thm}[N. Kovacevic]
\label{thm:semi-conjugacy}
Suppose $G<\Homeo(\R)$ is a non-abelian group with the
property that every nontrivial element of $G$ has at most one fixed point.
Then there is a continuous surjective semi-conjugacy from the action
of $G$ on $\R$ to a subgroup of the standard action of $\Aff(\R)$ on $\R.$
\end{thm}

We are indebted to E. Ghys for providing us with some history of these
two results. He attributes them both to V.V. Solodov, who apparently
never published a proof, but did at least announce a closely related
result in Theorem 3.21 of \cite{S}.  A published proof that a
non-abelian group with the property that every nontrivial element of
$G$ has at most one fixed point must be metabelian was provided by
T. Barbot in Theorem 2.8 of \cite{B}, at least in the case that each
fixed point is a topological attractor or a repellor.  Using deep
results on convergence groups N. Kovacevic published proofs of both
these results in \cite{K}.  She also considered the case that the
fixed points are topologically attractors or repellors.  It is likely
that both Barbot and Kovacevic understood how to deal with the minor
generalization to non-hyperbolic fixed points.

It is interesting to note that, prior to all of these results, J. Plante
proved a result in \cite{P} (and described below) which 
easily implies Theorem \ref{thm:semi-conjugacy} if one knows the group
$G$ is solvable.

In \S\ref{section:affineproof} we give a  proof of 
Theorem \ref{thm:affine} which takes H\"{o}lder's Theorem
as a starting point, and refines it by trying to reconstruct the
``translation subgroup'' of $G$ as the group of ``infinitesimals''.
In this way we prove $G$ is solvable.  Using the existence of a
quasi-invariant measure for $G$ (established by J. Plante \cite{P})
we then conclude the group is isomorphic to a subgroup of $\Aff(\R).$

Our main contribution to this study concerns the case when the action is
$C^2$, in which case we are able to show this semiconjugacy is actually a
topological conjugacy.  This gives a topological characterization of affine
groups among all groups of $C^2$ diffeomorphisms of $\R$.  E. Ghys
has subsequently informed us that he is also able to prove this
using techniques from Sacksteder's theory of pseudo-groups acting
on the line.

\begin{thm}[Conjugacy]
\label{thm:C2.conjugacy}
Suppose $G<\Diff^2(\R)$ is a non-abelian group with the property that
every nontrivial element of $G$ has at most one fixed point.  Then 
there exists $h\in \Homeo(\R)$ with $hGh^{-1}\in \Aff(\R)$.
\end{thm}

Theorem \ref{thm:C2.conjugacy} is not true for abelian groups, as can
be seen from an example of M. Hirsch \cite{H}.  In this paper Hirsch
gives an example of an analytic action of $\Z[1/2]$ on the line which
is semi-conjugate to the standard one, but which has an exceptional
minimal set and hence is not topologically conjugate to the standard
action.

\bigskip
\noindent
{\bf Actions on $\T^1$. }If $G$ is a group of homeomorphisms of the
circle $\T^1$, then the group $\widehat{G}$ of lifts to $\R$ of all
elements of $G$ is given by the exact sequence $$1\rightarrow
\Z \rightarrow \widehat{G}\rightarrow G\rightarrow 1$$ It follows easily from
H\"{o}lder's Theorem that any group of homeomorphisms of $\T^1$ which
acts freely must be abelian.  The following is a generalization of this
result, allowing elements to have a fixed point.

\begin{theorem} 
\label{theorem:circle} 
Let $G<\Homeo^+(\T^1)$ be torsion-free.  If every element of $G$ has at
most one fixed point, then $G$ is abelian.  Moreover, either $G$ acts
freely or $G$ has a global fixed point.
\end{theorem}  

Allowing elements of $G$ to have at most {\em two} fixed points
encompasses a much richer collection of examples, for example all
subgroups of $\PSL(2,\R)$.  

\section{Actions by elements with at most one fixed point}

In this section we consider a group $G$ acting on $\R$ with the
property that each nontrivial element has at most one fixed point.
We assume all elements preserve orientation of $\R.$

\subsection{Defining the order}

\begin{lemma}
Two elements of $G$ commute if and only if one of the following
conditions is satisfied:
\begin{enumerate}
\item The two elements generate a group acting freely on $\R.$
\item The two elements have a common fixed point.
\end{enumerate}
\end{lemma}

\begin{proof}
If the two elements, $g$ and $ h,$ commute then the fixed point set of
one is invariant under the other, so both have a fixed point or
neither does.  If both do it is a common fixed point.  If neither has
a fixed point then nothing in the group they generate has a fixed point
since the element with a fixed point would also commute with $g$ and $h.$

Conversely, if $g$ and $h$ generate a group acting freely, then they commute
by H\"older's theorem.  And if $g$ and $ h$ have a common fixed point $x_0,$
then they generate a group acting freely on the intervals $(-\infty, x_0)$
and $(x_0, \infty).$  By H\"older's theorem again we see they commute.
\end{proof}

\begin{lemma}\label{two_point}
If $g,h \in G$ and $g(x_0) = h(x_0)$ and $g(x_1) = h(x_1)$ with
$x_0 \ne x_1$ then $g = h$.  In other words, the graphs of distinct
elements of $G$ cross at most once.
\end{lemma}
\begin{proof}  The points $g(x_0)$ and $g(x_1)$ are two distinct
fixed points of $hg^{-1}$.  It follows that $hg^{-1}= id.$
\end{proof}

\begin{defn}
If $g,h \in G$ we say $g > h$ provided $g(x) > h(x)$ for all
sufficiently large $x,$ i.e. for all $x \in [x_0, \infty)$ for
some $x_0.$  We will call $g$ {\bf positive} provided $g>id$ and
{\bf negative} provided $id > g.$  
\end{defn}

Since the graphs of distinct elements of $G$ cross at most once, it
follows that for any distinct $g,h\in G$ either $g>h$ or $h>g$.  

\subsection{Commensurate elements}

\begin{defn}
We say that elements $h, g \in G$ are {\bf commensurate}
if there are $n,m\in \Z$ such that
$g^{-n} < h < g^n$ and $h^{-m} < g < h^m.$
\end{defn}

Note that taking inverses we see that $g^{-n} < h < g^n$ implies
$g^{-n} < h^{-1} < g^n$.  It follows that $h$ is commensurate with $g$
if and only if $h^{-1}$ is.  It is an easy exercise to see that
commensurability is an equivalence relation on $G$.

\begin{lemma}\label{KEY_LEMMA}
Suppose $g \in G$ is positive and $h \in G$ has a fixed point.
Then there exists $n \in \Z$ such that $h^n > g$ and $h^{-n} < g^{-1}.$
\end{lemma}

\begin{proof}
Replacing $h$ by $h^{-1}$ if necessary we may suppose
that $h$ is positive.  By conjugating $h$ we may assume that $0$ is
the fixed point of $h$ and that $h(x) = 2x$ for $x \ge 0$ and that for
$x < 0$ we either have $h(x) = 2x$ or $h(x) = x/2.$

Suppose first that for all $x$ we have $h(x) = 2x.$ We will show that
for large $n, \ h^n > g.$ Note that the graphs of $g$ and $h^n$ can
intersect in at most one point by Lemma \ref{two_point}.  Hence if
$g(0) = 0$ then the origin is the unique point of intersection of these
graphs.  Since $\lim_{n \to \infty} h^n(1) = \infty$ it follows that
for $n$  sufficienly large $h^n(1) > g(1)$ and hence $h^n > g.$
So we have the desired result if $g(0) = 0.$

If $g(0) \ge 0$ then $g(1) > g(0)$ and 
since $\lim_{n \to \infty} h^n(1) = \infty$ 
we may choose $N>0$ such
that $h^n(1) > g(1)$ for all $n > N.$
Since $h(0) = 0 \le g(0)$ there is an $x_n \in [0,1]$
with $h^n(x_n) = g(x_n)$.  
Note that the graphs of $g$ and $h^n$ can
intersect in at most one point by Lemma \ref{two_point}, so
the fact that $h^n(1) > g(1)$ implies that $h^n(x) > g(x)$ for
all $x > 1.$  Thus $h^n > g.$

Similarly if $g(0) < 0$ there is an $N$ such that
$h^N(-1) = -2^N < g(-1)$ and hence for any $n
> N$ there is an $x_n \in [-1,0]$ which is the unique point of
intersection of the graphs of  $h^n$ and of $g$.  
Since $h^n(0) > g(0)$ it follows that $h^n(x) > g(x)$ for all
$x > 0.$  So again $h^n > g.$

Consider now the case that $h(x) = x/2$ for $x < 0.$  
If $g(0) = 0$ then the argument above remains valid, so
we may assume the graph of $g$ intersects the $x$ and $y$ axes
in distinct points.

Letting $n$ go to infinity, the graphs of $h^n$ limit on the union of
positive $y$-axis and the negative $x$-axis.  Also letting $n$ go to
negative infinity these graphs limit on the union of positive $x$-axis
and the negative $y$-axis.  From this it is easy to see that there
exists an $n \in \Z$ such that the graph of $h^n$ intersects the graph
of $g$ in two points near the two points where the graph of $g$
intersects the axes.  But this is impossible by Lemma \ref{two_point}
\end{proof}

\begin{corollary}
\label{commensurate}
Let $g,h\in G$ be non-trivial.  If both $g$ and $h$ have fixed points
then then they are commensurate.
\end{corollary}

\begin{proof}
Recall that $h$ and $g$ are commensurate if and only if $h^{-1}$ and
$g$ are commensurate.  Hence replacing $h$ and $g$ by their inverses
if necessary we can assume both are positive.

By Lemma \ref{KEY_LEMMA}
there is an $n$ such that $h^n >g.$ So $h^{-n} < g^{-1}$ and we have
$h^{-n} < g^{-1} < g < h^n.$ The other inequality is similar.
\end{proof}

\subsection{Infinitesimals}

\begin{defn}
If $g \in G$ is positive we will define the set $\I(g)$ of 
elements of $G$ which are {\bf infinitesimal} with respect to $g$ by
\[
	\I(g) = \{ h \in G\ |\ g^{-1} < h^n < g \text{ for all } n \in \Z\}.
\]
If $g$ negative we define $\I(g) = \I(g^{-1}).$
\end{defn}

\begin{lemma}\label{independent}
If $g, g_1 \in G$ are commensurate
then $\I(g) = \I(g_1).$  Also if $h \in \I(g)$ and $h_1$ is commensurate
with $h$ then $h_1 \in \I(g).$
\end{lemma}
\begin{proof}
Without loss of generality we may assume both $g$ and $g_1$ are
positive. Suppose $k > 0$ and note that $g^k > g$ implies immediately
that  $\I(g) \subset \I(g^k)$.  But also $\I(g^k)\subset \I(g)$
since otherwise there would be an $h \in \I(g^k)$ and an $n>0$ 
with either $h^n > g$ or $h^n < g^{-1}$.  This would imply
$h^{nk} > g^k$ or $h^{nk} < g^{-k}$ contadicting the fact that
$h \in \I(g^k).$  We have shown that $\I(g) = \I(g^k).$

Since $g$ and $g_1$ are commensurate there is a $k>0$ such that $g^k >
g_1 > g_1^{-1} > g^{-k}.$ From this it is immediate that $\I(g_1)
\subset \I(g^k) = \I(g).$ Since $g$ and $g_1$ were arbitary positive
elements we conclude that $\I(g) = \I(g_1).$

To prove the second assertion of the lemma, suppose $h_1$ is
commensurate with $h$ and $h \in \I(g).$ We may assume both $h$ and
$h_1$ are positive.  Then there exists $k>0$ such that $h^{-k} <
h_1^{-1} < h_1 < h^k$.  It follows that $g^{-1} < h^{-nk} < h_1^{-n} <
h_1^n < h^{nk} < g$ for any $n >0$ and hence $h_1 \in \I(g).$
\end{proof}

We note that whenever $g$ and $g_1$ both have fixed points they are
commensurate by Corollary \ref{commensurate} and hence $\I(g) = \I(g_1)$ 
by Lemma \ref{independent}, so we can make the following definition.

\begin{defn}
We define $\I,$ the {\bf infinitesimals} of $G,$ to be $\I(g)$ if
$G$ contains a non-trivial element $g \in G$ which has a fixed point,
and by $\I = G$ if $G$ acts freely.
\end{defn}

\begin{proposition}\label{infinitesimal}
$\I$ is a normal abelian subgroup of $G$ which acts freely on $\R.$
\end{proposition}
\begin{proof}
First we note that if $G$ acts freely then the result is immediate
from H\"older's theorem and the definitions.  

Hence we may assume that $G$ contains a non-trivial positive element
$g$ with a fixed point.  Clearly we may assume there is a non-trivial
element $h \in \I(g).$

Since $h \in\I$ if and only if $h^{-1} \in\I$ we need only show
$\I$ is closed under multiplication in order to show it is a subgroup.
If $h_1, h_2 \in \I$ we must show $h_1h_2 \in \I.$ We may, without
loss of generality assume both are positive and that $h_2 > h_1$.
It follows that $h_1h_2 < h_2^2$ and hence that 
$(h_1h_2)^{-n} < id < (h_1h_2)^n < h_2^{2n} < g$
for all $n >0.$  So for all $n$ we have $(h_1h_2)^n < g$ and consequently
$g^{-1} < (h_1h_2)^{-n}.$  So $h_1h_2 \in \I.$

We note that Lemma \ref{KEY_LEMMA} implies
that any non-trivial element of $\I$ is fixed point free, so $\I$
acts freely and by H\"older's theorem it is abelian.

The fact that $\I$ is normal follows from the observation that
$h \in \I(g)$ and $a \in G$ implies $a^{-1} h a \in \I(a^{-1} g a).$
But as $g$ and $a^{-1} g a$ both have fixed points they are commensurate
by Corollary \ref{commensurate} and hence $\I(a^{-1} g a) = \I(g) =\I$
by Lemma \ref{independent}.
\end{proof}

\subsection{At most one fixed point implies solvable}

The following result of T. Barbot appears as Theorem
2.8 of \cite{B}.  It was apparently known 
to Solodov.

\begin{theorem}
\label{thm:solvable}
Let $G<\Homeo^+(\R)$.  Suppose that every nontrivial 
element of $G$ has at most one
fixed point.  Then $G$ is solvable with derived length at most two.
\end{theorem}

Since $\I$ is a normal Abelian subgroup it suffices to prove that $G/\I$
is Abelian since that implies that $[G,G] \subset \I.$ If $G$ acts
freely then $G = \I$ and the result is obvious, so we may assume there
exists $g \in G$ which has a fixed point and that $\I = \I(g).$

We define an order $\succ$  on $G/\I$ by $a\I \succ b\I$ if and only 
if $ab^{-1}$ is positive and not in $\I$.
To see this is well defined we must show that $ab^{-1} > id$ 
and $ab^{-1} \notin \I$ implies that
$(ah_1)(bh_2)^{-1} > id$ and $(ah_1)(bh_2)^{-1} \notin \I,$ 
for any $h_1, h_2 \in \I.$  

But $(ah_1)(bh_2)^{-1} = ab^{-1}h$ for some $h\in \I$ since $\I$ is
normal.  Clearly $ab^{-1} \notin \I$ implies that $ab^{-1}h \notin \I.$
So $(ah_1)(bh_2)^{-1} \notin \I.$ 

Also the fact that $ab^{-1} \notin \I$ implies $ab^{-1}$ is commensurate
with $g$.  Hence $\I = \I(g) = \I(ab^{-1})$ by Corollary \ref{commensurate}.
Thus $ab^{-1} > h^n$ for all $n$ and in particular for $n = -1$.  This
implies $ab^{-1}h > id.$  We have shown the order $\succ$ on $G/\I$ is well
defined.

We also observe that this order is both left and right invariant.  To
see this note that if $b_2\I \succ b_1\I$ then $ab_2\I \succ ab_1\I$
since $ab_1b_2^{-1}a^{-1}$ is a conjugate of $b_1b_2^{-1}$ and hence
positive and not in $\I.$ This proves left invariance; right
invariance is similar.

We next wish to show $G/\I$ is Archimedian with respect to the order
$\succ.$ Suppose $a\I$ and $b\I$ are positive elements of $G/\I$,
i.e. $a\I \succ \I$ and $b\I \succ \I$.  Then neither $a$ nor $b$ is
in $\I$ so they are both commensurate to $g$ and hence commensurate to
each other.  We conclude that there exists $k >0$ such that $a^k > b.$
Then $a^kb^{-1}$ is positive and if $a^kb^{-1} \in \I$ then
$a^{k+1}b^{-1}$ is positive and not in $\I.$ We have shown that in
either case $a^{k+1}\I \succ b\I.$ So $G/\I$ is Archimedian and hence
Abelian by H\"older's theorem.
\endproof

\subsection{Proof of Theorem \ref{theorem:onefixed}}
As before there exists a total order $<$ on $G$.  Since every element 
of $G$ has a fixed point, by Lemma \ref{KEY_LEMMA} we have $\I = \{id\}.$
Hence $G \cong G/\I$ is abelian by Theorem \ref{thm:solvable}.

\subsection{Proof of Theorem \ref{theorem:compact}}

As above we get a total order $<$ on $G$.  Given $g,h\in G$ their
graphs do not cross each other to the right of $b$.  Hence to prove
$g>h$ for $g,h\in G$ it is enough to prove $g(x_0)>h(x_0)$ for any
single point $x_0>b$.  By H\"{o}lder it is enough to check that, given
$h,g\in G$, $h^n>g$ for some $n\in \Z$.  By replacing $h$ with $h^{-1}$
if necessary, we may assume that $h$ is positive.  But then, for any 
fixed $x_0>b$, we have $h^n(x_0)\rightarrow \infty$, so given $g(x_0)$
simply pick $n$ big enough that $h^n(x_0)>g(x_0)$, which will imply 
$h^n>g$.

\subsection{Proof of Theorem \ref{theorem:circle}}

If $G$ does not act freely then there is an $h \in G$ with a fixed
point.  Call this fixed point $\infty$ and conjugate so $h(x) = x + 1$
on $(-\infty, \infty).$

Let $g$ be another element of $G$.  We want to show $g(\infty) = \infty$,
i.e. $\infty$ is a global fixed point.  If not, consider the graph of
$g$ on $\R = (-\infty, \infty).$  There is $x_0 \in \R$ such that
$g(x_0) = \infty$ and $y_0 \in \R$ such that $g(\infty) = y_0.$

This means the graph of $g$ has a vertical asymptote at $x = x_0$ and
a horizontal asymptote at $y = y_0$.  The graph is monotonically
increasing on $(-\infty, x_0)$ where it is entirely above the line 
$y = y_0$ and it is monotonically increasing on $(x_0, \infty)$ where it
is entirely below the line $y = y_0$.  That is, with respect to the
axes centered at $(x_0, y_0)$ the graph lies in the second and fourth
quadrants. 

The line $y = x - n$ is the graph of $h^{-n}$.  It intersects the
horizontal asymptote at $( y_0 + n, y_0)$ and the vertical 
asymptote at $( x_0 , x_0 - n)$.  Define $L_n(x) = h^{-n}(x) = x -n.$

For $n$ large we claim that the graph of this line must
intersect the graph of $g$ in at least two points.  To see this
choose a point $(x_1,y_1)$ below the graph of $g$ and with $x_1 > x_0$.
Then $(x_1,y_1)$ is above the graph of $y = L_n(x)$ for a sufficiently
large $n.$

Thus for small $\epsilon$ we have 
\begin{align*}
	g( x_0 + \epsilon) &< L( x_0 + \epsilon),\\
	L( x_1) &< g( x_1), \text{ and }\\
	g( y_0 + n ) &< L( y_0 + n).
\end{align*}
Hence the graphs intersect in at least one point in $(x_0,x_1)$ and
at least one point in $(x_1, y_0 +n)$.

Since we have found two points where $h^{-n}(x) = g(x)$ we have found
two fixed points of $h^ng$ and hence contradicted our hypothesis.  The
contradiction arose from assuming that $\infty$ was not a fixed point
of $g$.  Since $g$ was arbitrary we have shown that $\infty$ is a
global fixed point (under the assumption that one non-trivial element
has a fixed point).  Since $\infty$ is a global fixed point we can
consider $G$ as a free action on $\R = (-\infty, \infty)$.  By
H\"older's theorem it is abelian.

\section{Quasi-invariant measures}

A basic property of groups of homeomorphisms of $\R$ which act freely
(beyond being abelian) is that they have an invariant Borel measure.
This follows from the fact that the quotient space of $\R$ formed by
identifying orbits of a single element $h$ of $G$ is homeomorphic (or
diffeomorphic if $h$ is a diffeomorphism) to $\T^1$ and the quotient
group $G/<h>$ acts on this $\T^1$.  This action has an invariant Borel
probability measure which may be lifted to $\R$ to give an invariant
Borel measure $\mu$ for $G$ which is finite on compact sets and 
infinite on all of $\R.$

The following result of J. Plante (Theorem 4.7 of \cite{P}) asserts
that this can be extended to a quasi-invariant measure for certain
solvable groups.

\begin{thm}[J. Plante]
\label{thm:quasi.invariant}
Suppose $G<\Homeo(\R)$ is solvable, and its nontrivial elements have
isolated fixed points.  Then there exists a $G$-quasi-invariant measure
$\mu$ which is finite on compact sets.
\end{thm}

A measure is $G$ {\em quasi-invariant} if
for each $g \in G$ the measure $g_*(\mu)$ is
equal to $A(g)\mu$ for some positive real number $A(g)$.  The
function $A: G \to \R^+$ is a homomorphism to the multiplicative
group $\R^+$.  The measure $\mu$ is invariant under the subgroup
which is the kernel of $A.$

If one has an invariant measure $\mu$  for a subgroup of $\Homeo(\R)$ it is
useful to consider the translation number which was discussed by 
J.~Plante in \cite{P}.  In our context this is just the rotation number for
circle homeomorphisms obtained by quotienting $\R$ by an element of 
$H$. 

Suppose $H$ is a subgroup of $\Homeo(\R)$ which preserves a Borel measure
$\mu$ that is finite on compact sets.  Fix a point $x \in \R$ and for each $f \in G$
define 
\begin{equation*}
\tau_\mu(f) =
\begin{cases}
\mu([x,f(x)))	&\text{if $x < f(x)$} \\ 
0	&\text{if $x =f(x)$ }\\
-\mu([f(x),x))	&\text{if $x >f(x)$} 
\end{cases}
\end{equation*}

The function $\tau_{\mu} : G \to \R$ is called the {\em translation number}.
The following properties observed by J. Plante in \cite{P} are easy
to verify.

\begin{prop}\label{prop:translation.number}
The translation number $\tau : G \to \R$ is independent
of the choice of $x \in \R$ used in its definition.
It is a homomorphism from $G$ to the additive group $\R.$
For any $f \in G$ the set $\Fix(f) \ne \emptyset$
if and only if $\tau_\mu(f) = 0.$
\end{prop}

If $G$ is a subgroup of $\Homeo(\R)$ with the property that
every non-trivial element has at most one fixed point then we
can combine the homomorphisms $\tau$ and $A$ from \cite{P}
to obtain a homomorphisms to the affine group $\Aff(\R).$

\begin{lem}
\label{lem:kernel.A}
Suppose $G<\Homeo(\R)$ is a non-abelian group with the
property that every nontrivial element has at most one fixed point.
If $\mu$ is a $G$-quasi-invariant measure which is finite on compact
sets, then a nontrivial element 
$g$ has exactly one fixed point if and only if $A(g) \ne 1$.
\end{lem}

\begin{proof}
If $ A(g) = 1$ then $\mu$ is $g$-invariant and $\tau_{\mu}(g)$ is well
defined.  If $\tau_{\mu}(g)\ne 0$ then $g$ has no fixed points by
Proposition \ref{prop:translation.number}.  If $\tau_\mu(g) = 0$ then
for each $x \in {\R},\ \mu([x, g^n(x)))$ is bounded as a function of
$n$.  This implies $g$ has infinitely many fixed points and hence must
be the trivial element.  

Suppose now that $A(g) = \alpha \ne 1$.  Replacing $g$ by $g^{-1}$ if
necessary, we may assume $0 < \alpha < 1$.  We will show that the
assumption that $g$ has no fixed points leads to a contradiction.  Let
$x \in \R$ and assume without loss of generality that $g(x) < x$ for all
$x \in \R$, so $\lim_{n \to \infty}g^n(x) = -\infty.$ Then
\[
\mu([g^n(x), x]) \le 
\sum_{i=0}^{n-1} \mu([g^{i+1}(x),g^i(x)]) 
\le \sum_{i=0}^{\infty}\alpha^i \mu([g(x),x]) =  
\frac{\mu([g(x), x])}{1-\alpha}.
\]

Hence for every $x, \ \mu((-\infty, x])$ is finite.  But this is not
possible, as follows.  Since $G$ is non-abelian, there is an 
element $h$ in the kernel of $A$.  The element $h$ has the
property that $\lim_{n \to \infty}h^n(x) = -\infty$ and $h$ preserves
the measure $\mu$.  Thus
\[
\mu((-\infty, x]) = \sum_{i=0}^{\infty} \mu((h^{i+1}(x),h^i(x)]) 
= \lim_{n \to \infty} n \mu((h(x),x]),
\]
which imples that for every $x$ that $\mu((-\infty, x])$ is either 
$0$ or infinite, a contradiction.

Thus we conclude that $g$ has a fixed point.  Since $A(g) \ne 1$ it is not the
identity, so it has exactly one fixed point.
\end{proof}

\subsection{Proof of Theorem \ref{thm:affine}}
\label{section:affineproof}

We define a function $\nu : \R^2 \to \R$ by
\begin{equation*}
\nu( x_0,x_1)  =
\begin{cases}
\mu([x_0,x_1))	&\text{if $x_0 <x_1$} \\
0&\text{if $x_0 = x_1$} \\ 
-\mu([x_0,x_1))	&\text{if  $x_0 > x_1$}
\end{cases}
\end{equation*}
Then define
$\phi: G \to \Aff(\R)$ by $\phi(g)(x) = ax + b$, where
$a = A(g)$ and $b = \nu(0,g(0)).$

If $g_1,g_2 \in G$ then 
\[
\phi( g_1 g_2)  = A(g_1) A( g_2) x + \nu(0, g_1g_2(0)).
\]
But 
\begin{align*}
\nu(0, g_1g_2(0)) &= \nu(0, g_1(0)) + \nu(g_1(0), g_1g_2(0))\\
&= \nu(0, g_1(0)) + A(g_1) \nu(0, g_2(0)).
\end{align*}
So 
\begin{align*}
\phi( g_1 g_2)  &= A(g_1)( A( g_2) x + \nu(0, g_2(0)) + \nu(0, g_1(0))\\
&=\phi(g_1) \phi(g_s).
\end{align*}

To see that $\phi$ is injective, observe that if $g$ is in the kernel
of $\phi$ then $A(g) = 1$ and $\phi(g)(x) = x + \nu(0, g(0)) = x +
\tau_{\mu}(g)$, so that $\tau_{\mu}(g) = 0.$ Now if $G$ is abelian then 
the theorem is clearly true, so suppose $G$ is non-abelian.  From Proposition
\ref{prop:translation.number} we conclude that $g$ has a fixed point
and from Lemma \ref{lem:kernel.A} that it has more than one.  Hence $g
= id$ and $\phi$ is injective.

\subsection{Proof of Theorem \ref{thm:semi-conjugacy}}

$G$ is solvable by Theorem \ref{theorem:onefixed}.  Hence by Theorem
\ref{thm:quasi.invariant} there is a $G$-quasi-invariant measure $\mu$.

Define $\theta(x) = \nu(0,x).$
It was observed by J. Plante (see the remark after 4.6 of \cite{P})
that in this case $\theta$ is a semi-conjugacy, as the following
calculation shows.
\begin{align*}
\theta(g(x)) &= \nu(0,g(x))\\
&= \nu(0,g(0)) + \nu(g(0),g(x))\\
&=\nu(0,g(0)) + A(g) \nu(0,x)\\
&= A(g)((\theta(x)) + \nu(0, g(0))\\
&=\phi(g)(\theta(x)).
\end{align*}
Hence $\theta(g(x)) = \phi(g)(\theta(x))$ for all $x \in \R$ and $g \in G.$

\section{$C^2$ actions}

In this section we prove that the semi-conjugacy in Theorem
\ref{thm:semi-conjugacy} is actually a topological conjugacy in the
case that $G$ is a subgroup of $\Diff^2(\R)$.p

\begin{defn}
\label{defn:distortion}
If $f: \R \to \R$ is a diffeomorphism and $J$ is an interval in $\R$ then
we define the {\em distortion} of $f$ on $J$ by
\[
\Dist(f,J) = \sup_{x,y \in J}\log\frac{ |Df(x)|}{ |Df(y)|}.
\]
\end{defn}

The following lemma is well known; a proof can be found in Chapter I, \S2 of 
\cite{dMvS}.  We use the notation $|J|$ to denote the length of an interval
$J.$

\begin{lem}\label{lem:distortion}
If $f: \R \to \R$ is a $C^2$ diffeomorphism of $\R$ and the function
$\log |Df(x)|$ is Lipschitz with Lipschitz constant $C$, 
then for all $n>0,$
\[
\Dist(f^n, J) \le C \sum_{i=0}^{n-1} | f^i(J)|.
\]
\end{lem}

We will also need the following lemma.  This result and the proof we
give are extracted from a part of the proof a result of A. Schwartz as
presented by de Melo and van Strien in Theorem 2.2 of Chapter I, \S2
from \cite{dMvS}.

\begin{lem}\label{lem:wandering.bound}
Suppose $g: \R \to \R$ is a $C^2$ diffeomorphism and
$\log |Dg(x)|$ is Lipschitz with a uniform Lipschitz constant $C,$  and suppose $J \subset \R$
is a closed interval with the property that
\[
\sum_{i=0}^{\infty} |g^i(J)| \le 1.
\]
Then there is a $\delta >0$ such that any closed interval $L$ which contains $J$
and has length $|L| < (1 + \delta)|J|$ has the property that
\[
\sum_{i=0}^{\infty} |g^i(L)| \le 2.
\]
\end{lem}

\begin{proof}
Choose  $\delta$ so $0 < \delta < \min\{ |J|, \exp(-2C) \}$.
We will prove by induction that $|g^n(L)| \le 2|g^n(J)|$ for any $n \ge 0.$
This is clear for $n= 0.$  Assume inductively that
 $|g^k(L)| \le 2|g^k(J)|$ for any $k <n.$

We know by Lemma \ref{lem:distortion} and the induction hypothesis that 
\begin{equation}
\Dist(g^n, L) \le C \sum_{i=0}^{n-1} | g^i(L)| \le 2 C \sum_{i=0}^{n-1}
 | g^i(J)| \le 2C. \label{dist_eq}
\end{equation}
Also by the mean value theorem there is a point $x_n \in J$ such that
$|Dg^n(x_n)| = |g^n(J)|/|J|.$  Therefore by equation (\ref{dist_eq})
for every $y \in L$
\[
|Dg^n(y)| \le \exp(2C)|Dg^n(x_n)| \le \exp(2C) \frac{|g^n(J)|}{|J|}.
\]

Using the mean value theorem again and the fact that  $(|L| - |J|)/|J| < \delta$, 
we see
\begin{align*}
|g^n(L)| &\le |g^n(J)| + \sup_{y \in L}|Dg^n(y)| (|L| - |J|) \\
&\le |g^n(J)| + \exp(2C) \frac{|g^n(J)|}{|J|}(|L| - |J|) \\
&\le (1 + \delta \exp(2C)) |g^n(J)|\\
& \le 2 |g^n(J)|.
\end{align*}

This completes the induction step, so $|g^n(L)| \le 2 |g^n(J)|$ for
all $n \ge 0.$ The fact that $\sum_{i=0}^{\infty} |g^i(L)| \le 2$ then
follows from the fact that $\sum_{i=0}^{\infty} |g^i(J)| \le 1.$
\end{proof}

\subsection{Proof of Theorem \ref{thm:C2.conjugacy}}

{\bf Theorem \ref{thm:C2.conjugacy}} {\em Suppose $G<\Diff^2(\R)$ is a
non-abelian group with the property that every nontrivial element of $G$
has at most one fixed point.  Then the function $\theta: \R \to \R$
defined by $\theta(x) = \nu(0,x)$ is a topological conjugacy from the
action of $G$ on $\R$ to the standard action of $\Aff(\R)$ on $\R.$}

\begin{proof}
In Theorem \ref{thm:semi-conjugacy} we showed the function $\theta: \R
\to \R$ defined by $\theta(x) = \nu(0,x)$ is a surjective
semi-conjugacy.  It suffices to prove that $\theta$ is one-to-one.
Note that if $\mu$ is the $G$-quasi-invariant measure provided by
Theorem \ref{thm:quasi.invariant}, then $\theta_{\ast}(\mu)$ is
Lebesgue measure on $\R$, the unique quasi-invariant measure for the
standard action of $\Aff(\R)$ on $\R.$ The function $\theta$ is
monotonic so $\theta^{-1}(y)$ consists of either a single point or a
closed interval.  If it is an interval $[x_0,x_1]$ then
$\mu([x_0,x_1]) = 0$.  The endpoints $x_0,x_1$ are in the support of
$\mu$ but the interior $(x_0,x_1)$ is disjoint from the support of
$\mu.$

We note that any interval $[x_0,x_1] = \theta^{-1}(y)$ must be
wandering.  That is, we must have $g([x_0,x_1]) \cap [x_0,x_1] =
\emptyset$ for all non-trivial $g \in G.$ This is because
if $g([x_0,x_1]) \cap [x_0,x_1] \ne \emptyset$ the fact that $\theta$
is a semi-conjugacy implies $g([x_0,x_1]) = [x_0,x_1]$ so $g$ would
have two fixed points $x_0$ and $x_1$.

Thus it suffices to show that the action of $G$ has no wandering
intervals of the form $\theta^{-1}(y)$.  Of course to do this we need
only do it for a subgroup of $G$.  In the proof of Theorem
\ref{thm:affine} we showed that the function $\phi: G \to \Aff(\R)$
defined by $\phi(g)(x) = ax + b$, where $a = A(g)$ and $b =
\nu(0,g(0)),$ is an injective homomorphism.  In particular if $A(h) =
1$ then $\phi(h)(x) = x + \tau_{\mu}(h).$ Since we are assuming $G$ is
not abelian, $A: G \to \R^+$ is not trivial and neither is its kernel.
Hence we will focus on a subgroup $G_0$ generated by two non-trivial
elements, $g$ with $A(g) > 1,$ and $h$ with $A(h) = 1$ and
$\tau_{\mu}(h) \ne 0$.  Replacing $h$ by $h^{-1}$ if necessary, we may
assume $h(x) > x$ for all $x$.

Consider the quotient space of $\R$ by the action
of $h$.  There is a $C^2$ diffeomorphism from this quotient to $\T^1$.
We lift this diffeomorphism to a $C^2$ diffeomorphism of $\R$ and use
it to conjugate our $G$ action to a $C^2$ action on $\R$ for which
$h(y) = y + 1.$ We know $\phi(h) = x + c$ for some $c >0$ but
composing $\phi$ with an inner automorphism of $\Aff(\R)$ (i.e. an
affine change of co-ordinates) we may assume $\phi(h) = x + 1$ and
$\phi(g) = \lambda x$ for some $\lambda > 1.$

We first consider the case that $\lambda$ is irrational.  Let $f =
ghg^{-1}$.  Since $h$ is in the kernel of $A$, a normal abelian
subgroup, we have that $f$ commutes with $h$.  We also have 
that $\phi(f) = x + \lambda$.  
Hence $f$ is the lift of a $C^2$ diffeomorphism of the circle $\T^1$
with irrational rotation number.  By Denjoy's theorem this circle
diffeomorphism has a dense orbit from which it follows that if $H$ is
the group generated by $h$ and $f$ then the set $H(x)$ is dense in
$\R$ for any $x$.  Hence $H$ has no wandering intervals at all.

Thus we may assume $A(g) = p/q$ is rational.  Note that $A(g) = p/q$
implies $gh^q = h^pg$ or $g( y + q) = g(y) + p.$ Differentiating we
see that for all $y \in \R$ we have $Dg(y + q) = Dg(y)$ and $D^2g(y +
q) =D^2g(y)$.  In particular $Dg(y)$ and the derivative of $\log
Dg(y)$ are periodic and hence uniformly bounded.

We now want to show by contradiction that $G_0$ has no wandering
intervals of the form $\theta^{-1}(y)$.  Suppose $J = [x_0,x_1] =
\theta^{-1}(y)$ is a wandering interval.  Then for every non-trivial
$g \in G_0, \ \phi(g)(y) \ne y$ because otherwise $g(\theta^{-1}(y)) =
\theta^{-1}(y)$ and the endpoints of this $g$ invariant interval would
be two fixed points for $g.$

We can conclude, in particular, that the elements $y_n = \phi(g)^n(y)
\mod(1) = \phi(g^n)(y) \mod(1)$ are all distinct, where $y \mod(1)$
denotes $y$ minus the greateset integer in $y$.  This is because
otherwise $g(J)$ would equal $J$ for some $g \in G.$
It follows that the intervals 
$J_n = \theta^{-1}(y_n)$ are all pairwise disjoint.

Since $\theta$ is a semiconjugcay and is $h$-equivariant,
$\theta(x+n)=\theta(x)+n$ for all $n\in \Z$.  As $\theta$ is monotonic
and $\theta^{-1}(n)=n$ for all $n \in \Z$, we have
$\theta^{-1}([n,n+1])=[n,n+1]$.  In particular
$\theta^{-1}([0,1])=[0,1],$ so since $J_k = \theta^{-1}(y_k)$ and $y_k
\in [0,1]$ we have $J_k \subset [0,1]$ for all $k >0.$ Recall that
$|J_k|$ denotes the length of $J_k$ and observe that $J_k = g^k(J) -
n_k$ for some integer $n_k$, so $|J_k| = |g^k(J)|.$ Hence
\[
\sum_{k=0}^\infty |g^k(J)|  = \sum_{k=0}^\infty |J_k| \le 1,
\]
as the intervals $J_k$ are pairwise disjoint and all in $[0,1].$

Let $\delta > 0$ be the value provided by Lemma
\ref{lem:wandering.bound} and choose an interval $L$ which contains
$J$ in its interior and has length $|L| < (1+\delta)|J|.$ Define $L' =
\theta(L)$ and note that it is a non-trivial closed interval with $y =
\theta(J)$ in its interior.  We may shrink $L'$ slightly so
that $\theta^{-1}$ of each endpoint is a single point (because
$\theta^{-1}(x)$ is a non-trivial interval for at most
countably many values of $x.$) We redefine $L$ to be $\theta^{-1}$ of
this shrunken $L'$ and note that $L = \theta^{-1} (L') = \theta^{-1}
(\theta(L)).$ It follows that $g^n(L) = \theta^{-1}( \lambda^n L')$
for all $n>0.$

We observe that from Lemma \ref{lem:wandering.bound} that
\begin{equation}
\sum_{k=0}^\infty |g^k(L)|  \le 2. \label{eqn1}
\end{equation}

But $\theta( g^k(L)) = \lambda^k L'$ and for $k$ sufficiently
large $\lambda^k |L'| > 2$ which implies there is an interval
$[m,m+1] \subset \lambda^k L',$ with $m \in \Z.$ This means
that 
\[
[m,m+1] = \theta^{-1}([m,m+1]) \subset \theta^{-1}(\lambda^k L') =g^k(L). 
\]
Hence $|g^k(L)| \ge 1$ for $k$ sufficiently large.  This clearly
contradicts equation (\ref{eqn1}) above.  We conclude that the
semi-conjugacy $\theta$ is one-to-one and hence a homeomorphism.
\end{proof}

\bigskip
\noindent
Benson Farb:\\
Dept. of Mathematics, University of Chicago\\
5734 University Ave.\\
Chicago, Il 60637\\
E-mail: farb@math.uchicago.edu
\medskip

\noindent
John Franks:\\
Dept. of Mathematics, Northwestern University\\
Evanston, IL 60208\\
E-mail: john@math.northwestern.edu


\begin{thebibliography}{ABCDEF}


\bibitem[B]{B}
T. Barbot, Caracterisation des flots d'Anosov en dimension 3 par leurs
feuilletages faibles, {\em Ergodic Theory Dynam.
Systems} {\bf 15} (1995), no. 2, 247--270.

\bibitem[dMvS]{dMvS}
Welington de Melo and Sebastien van Strien
{\em One-dimensional Dynamics} 
Springer Verlag, Berlin, (1991)

\bibitem[FF1]{FF1}
B. Farb and J. Franks,
Group actions on one-manifolds, I: Nonlinear subgroups, 
in preparation.

\bibitem[FS]{FS}
B. Farb and P. Shalen, Groups of real-analytic diffeomorphisms of the
circle, to appear in {\em Ergodic Theory and Dynam. Syst.}

\bibitem[Gh]{Gh}
E. Ghys, Groups acting on the circle, IMCA, Lima, June 1999.

\bibitem[H]{H}
M. Hirsch, A stable analytic foliation with only exceptional minimal
set, in {\em Lecture Notes in Math.}, Vol. 468. Springer-Verlag, 1975.

\bibitem[Ho]{Ho}
O. H\"{o}lder, Die Axiome der Quantit\"{a}t und die Lehre vom
Mass. Ber. Verh. Sachs. Ges. Wiss. Leipzig, Math. Phys. C1. 53,
1-64 (1901).

\bibitem[K]{K}
N. Kovacevic,  M\"{o}bius-like groups of homeomorphisms of the circle.
{\em Trans. Amer. Math. Soc.} {\bf 351} (1999), no. 12, 4791--4822. 

\bibitem[P]{P}
J. Plante, Solvable Groups acting on the line,
{\em Trans. Amer. Math. Soc.} {\bf 278}, (1983) 401--414.

\bibitem[P2]{P2}
J. Plante, Subgroups of continuous groups acting differentiably on the
half-line, {\em Ann. Inst. Fourier, Grenoble} {\bf 34} 1 (1984), 47-56.

\bibitem[S]{S}
V. V. Solodov, Topological problems in the theory of dynamical systems.
(Russian) {\em Uspekhi Mat. Nauk} 46 (1991), no. 4(280),93--114, 192 
translation in {\em Russian Math. Surveys} 46 (1991), no. 4, 107--134

\end{thebibliography}
\end{document}